\RequirePackage{fix-cm}
\documentclass{svjour3}                     
\smartqed  
\usepackage{graphicx}
 \usepackage{mathptmx}      
\usepackage{amsmath,amsxtra,amssymb,latexsym,amscd}

\begin{document}

\title{ Extragradient algorithms for equilibrium problems and symmetric generalized hybrid mappings}
\author{Bui Van Dinh        \and Do Sang Kim    }
\institute{Bui Van Dinh  \at
              Faculty of Information Technology, Le Quy Don Technical University, Hanoi, Vietnam \\
              \email{vandinhb@gmail.com}   \\        
             \emph{Present address:} Department of Applied Mathematics, Pukyong National University, Busan, 608-737, Korea  
           \and
          Do Sang Kim \at
              Department of Applied Mathematics, Pukyong National University, Busan, 608-737, Korea\\
\email{dskim@pknu.ac.kr} }
\date{Received: date / Accepted: date}
\maketitle

\begin{abstract}
In this paper, we propose  new algorithms for finding a common point of the solution set of a pseudomonotone equilibrium problem and the set of fixed points of a symmetric generalized hybrid mapping in a real Hilbert space. The convergence of the iterates generated by each method is obtained under assumptions that the fixed point mapping is quasi-nonexpansive and demiclosed at $0$, and the bifunction associated with the equilibrium problem is weakly continuous. The bifunction is assumed to be satisfying a Lipschitz-type condition when the basic iteration comes from the extragradient method. It becomes unnecessary when an Armijo back tracking linesearch is incorporated in the extragradient method.

 \subclass{$47$H$06$ \and  $47$H$09$ \and  $47$H$10$ \and  $47$J$05$ \and  $47$J$25$}
\end{abstract}

\section{Introduction}

 Let $\mathbb{H}$ be a real Hilbert space with the inner product $\langle \cdot , \cdot\rangle$ and induced norm $\| \cdot \|$. By `$\to$'  and `$\rightharpoonup$' we denote the strong convergence and the weak convergence in $\mathbb{H}$, respectively. Let $C$ be a nonempty closed convex subset of $\mathbb{H}$ and $f : C \times C \to \mathbb{R}$ be a bifunction satisfying $f(x, x) = 0$ for every $x \in C$. Such a bifunction is called an equilibrium bifunction. The equilibrium problem,  in the sense of Blum, Muu and Oettli \cite{BO,MO} (shortly EP($C, f$)), is to find $ x^* \in C $ such that
\begin{equation}\label{1.1}
  f(x^*, y) \geq 0, \  \forall y \in C.
\end{equation}
By $Sol(C, f)$, we denote the solution set of EP($C, f$).
Although problem EP($C, f$) has a simple formulation, it includes, as special cases, many important problems in applied mathematics: variational inequality problem, optimization problem, fixed point problem, saddle point problem, Nash equilibrium problem in noncooperative game, and others; see, for example, \cite{BCPP,BO,MO}, and the references quoted therein.

Let us denote the set of fixed points of a mapping $T : C \to C$ by $Fix(T)$; that is, $Fix(T) = \{x \in C: Tx = x\}$. Recall that  $T$ is said to be nonexpansive if for all $x, \  y \in C$, \  $\|Tx - Ty\| \leq \|x - y\|$. If $Fix(T)$ is nonempty and $\|Tx - p \| \leq \| x - p \|, \ \forall x \in C, \  p \in Fix(T)$, then $T$ is called quasi-nonexpansive. It is well-known that $Fix(T)$ is closed and convex when $T$ is quasi-nonexpansive \cite{IT}.

A mapping $T$ is said to be  pseudocontractive if for  all $x, \  y \in C$ and $\tau > 0$,
 $$\|x - y\| \leq \|(1+\tau)(x - y) - \tau(Tx - Ty)\|.$$

To find a fixed point of a Lipschitzian pseudocontractive map, Ishikawa \cite{Ish}, in 1974, proposed to use the following iteration procedure
\begin{equation}\label{1.2}
\begin{cases}
x^0 \in C,\\
y^k =  \alpha_k x^k + (1 - \alpha_k)Tx^k, \\
x^{k+1} = \beta_k x^k + (1 - \beta_k)Ty^k
\end{cases}
\end{equation}
where $0 \leq \alpha_k \leq \beta_k \leq 1$ for all $k$ and proved that if $\lim_{k \to \infty} \beta_k = 1$, $\sum_{k=1}^{\infty}(1 - \alpha_k)(1 - \beta_k) = \infty$, then $\{x^k\}$ generated by (\ref{1.2}) converges weakly to a fixed point of mapping $T$ (see \cite{GL,Ish}).

In 2006, Yanes and Xu \cite{YX} introduced the following  by combining Ishikawa iteration process with hybrid projection method \cite{NT} for a nonexpansive mapping $T$.
\begin{equation}\label{1.3}
\begin{cases}
x^0 \in C,\\
y^k =  \alpha_k x^k + (1 - \alpha_k)Tx^k, \\
z^{k} = \beta_k x^k + (1 - \beta_k)Ty^k,\\
C_k = \{x \in C: \|x - z^k \|^2  \leq \|x - x^k \|^2 + (1 - \alpha_k)(\|y^k \|^2 - \| x^k \|^2  + 2\langle x^k - y^k, x \rangle)  \},\\
Q_k = \{x \in C: \langle x - x^k , x^0 - x^k \rangle \leq 0 \}, \\
x^{k+1} =  P_{C_k \cap Q_k} x^0,
\end{cases}
\end{equation}
where $\{\alpha_k\}$  and $\{\beta_k\}$ are sequences in [0, 1]. They proved that if $\lim_{k \to \infty} \alpha_k = 1 $ and $\beta_k \leq \bar{\beta}$  for some $\bar{\beta} \in [0, 1)$, then $\{x^k\}$ generated by (\ref{1.3}) converges strongly to $P_{Fix(T)} (x^0).$

In recent years, many researchers studied the problem of finding a common element of the set of solutions of an equilibrium problem and the set of fixed points of a nonexpansive or demicontractive mapping; see, for instance,  \cite{AM,CAY,Mai2,PK,VSH1} and the references therein.
Remember that a mapping $T : C \to \mathbb{H}$ is called symmetric generalized hybrid \cite{HST,KT,TWY} if there exist $\alpha, \beta, \gamma, \delta \in \mathbb{R}$ such that
$$ \alpha\|Tx - Ty\|^2 + \beta\big (\|x - Ty\|^2 + \|y - Tx \|^2 \big) + \gamma\|x - y\|^2 + \delta\big (\|x - Tx\|^2 + \|y - Ty \|^2 \big) \leq 0, \forall x, y \in C.$$ Such a mapping is called an $(\alpha, \beta, \gamma, \delta)$-symmetric generalized hybrid mapping.

For obtaining a common element of the set of solutions of EP($C, f$) and fixed points of a symmetric generalized hybrid mapping $T$, Moradlou and Alizadeh \cite{MA} proposed to combine Ishikawa iterative scheme with the hybrid projection method \cite{NT,TT,TTK}. More precisely, the iterates $x^k, \ y^k, \  u^k,\  z^k$ are calculated as follows:
\begin{equation}\label{1.4}
\begin{cases}
x^0 \in C,\\
u^k \in C \text{ such that } f(u^k, y) + \frac{1}{r_k}\langle y - u^k, u^k - x^k \rangle \geq 0, \forall y \in C,\\
y^k =  \alpha_k x^k + (1 - \alpha_k)Tx^k, \\
z^{k} = \beta_k y^k + (1 - \beta_k)Tu^k,\\
C_k = \{x \in C: \|x - z^k \|  \leq \|x - x^k \|  \},\\
Q_k = \{x \in C: \langle x - x^k , x^0 - x^k \rangle \leq 0 \}, \\
x^{k+1} =  P_{C_k \cap Q_k} x^0.
\end{cases}
\end{equation}
The authors showed that if $Sol(C, f) \cap Fix(T) \neq \emptyset$, $(\alpha, \  \beta, \  \gamma, \  \delta)$-symmetric generalized hybrid mapping $T$ satisfying $(1) \ \alpha + 2\beta + \gamma \geq 0$, $(2) \ \alpha + \beta > 0 $, and $(3) \ \delta \geq 0$, then under certain appropriate conditions imposed on $\{\alpha_k\}$, $\{\beta_k\}$, the sequence $\{x^k\}$  converges strongly to  $x^* = P_{Sol(C, f) \cap Fix(T)}(x^0)$ provided that $f$ is monotone on $C.$ \\
Note that mapping $T$ satisfies the conditions $(1), \ (2)$, and $(3)$, then $T$ is quasi-nonexpansive and demiclosed at $0$.

In this paper, we modify Moradlou and Alizadeh's iteration process for finding a common element of the set of solutions of an equilibrium problem and the set of fixed points of a generalized hybrid mapping  in a real Hilbert space  in which the bifunction $f$ is pseudomonotone on $C$ with respect to $Sol(C, f)$. More precisely, we propose to use the extragradient algorithm \cite{Kor} for solving the equilibrium problem (see also \cite{DM,DHM,FP,Kon,MQH} for more detail extragradient algorithms). One advantage of our algorithm is that it could be applied for the pseudomonotone equilibrium problem case and each iteration we only have to solve two strongly convex optimization problems instead of a regularized equilibrium as in Moradou and Alizaded's method.

The paper is organized as follows. The next section contains some preliminaries on the metric projection, equilibrium problems and symmetric generalized hybrid mappings. An extragradient algorithm and its convergence is presented in the third section. The last section is devoted to presentation of  an extragradient algorithm with linesearch and its convergence.


\section{Preliminaries}

In the rest of this paper, by $P_C$ we denote the metric projection operator on $C$, that is
 $$P_C(x) \in C: \ \|x - P_C(x) \| \leq \|y - x\|, \ \forall y \in C,$$

 The following well known results on the projection operator onto a closed convex
 set will be used in the sequel.

\begin{lemma}\label{Lem2.1} Suppose that $C$ is a nonempty closed convex subset in $\mathbb{H}$. Then
\begin{itemize}
\item[(a)] $P_C(x)$ is singleton and well defined for every $x$;

\item[(b)] $z = P_C(x)$ if and only if  $ \langle x - z  , y - z  \rangle \leq 0, \forall y \in C;$

\item[(c)] $\|P_C(x) - P_C(y)\|^2  \leq \|x-y\|^2 - \|P_C(x) - x + y - P_C(y)\|^2 ,  \ \forall x, y \in C$.
\end{itemize}
 \end{lemma}

\begin{lemma}\label{Lem2.4}\cite{YX}
Let $C$ be a nonempty closed convex subset of $\mathbb{H}$. Let  $\{x^k\}$ be a sequence in $\mathbb{H}$ and $u \in \mathbb{H}$. If any weak limit point of $\{x^k\}$ belongs to $C$ and $$ \|x^k - u \| \leq \|u - P_C(u)\|, \ \forall k.$$
Then $x^k \to P_C(u)$.
\end{lemma}

\begin{definition}\label{Def2.2} A  bifunction  $\varphi : C \times C \to \mathbb{R}$ is said to be
jointly weakly continuous on $C \times C$  if for all $x, y \in C$ and $\{x^k\}$,  $ \{y^k\} $ are two sequences in $C$ converging weakly to $x$ and $y$ respectively, then $\varphi(x^k, y^k)$ converges to $\varphi(x, y)$.

\end{definition}
In the sequel, we need the following blanket assumptions

\begin{itemize}
\item[$(A_1)$] $f$ is jointly weakly continuous on $C\times C$;
\item[$(A_2)$]  $f(x, \cdot)$ is convex, lower semicontinuous, and subdifferentiable on $C$, for all $x\in C$;
\item[$(A_3)$] $f$ is pseudomonotone on $C$ with respect to $Sol(C,f)$, i.e., $f(x, x^*)\leq 0$ for all $x \in C$, $x^* \in Sol(C, f)$;
\item[$(A_4)$] $f$ is Lipschitz-type continuous on $C$ with constants $L_1 >0$ and $L_2 >0$, i.e.,
$$f(x, y)+f(y, z)\geq f(x, z) - L_1 \Vert x-y \Vert^2 - L_2 \Vert y-z \Vert^2, \; \forall x, y, z\in C;$$
\item[$(A_5)$] $T$ is an ($\alpha, \beta, \gamma, \delta$)-symmetric generalized hybrid self-mapping of $C$ such that $(1) \ \alpha + 2\beta + \gamma \geq 0$,  $(2) \ \alpha + \beta > 0$,  $(3) \ \delta \geq 0 $, and $Fix(T)$ is nonempty.

    \end{itemize}

For each   $z, \  x \in C$, by $\partial_2 f(z , x)$ we denote the subgradient of the convex
  function $f(z, .)$ at $x$, i.e.,
$$ \partial_2 f(z , x) := \{ w \in \mathbb{H} : f(z, y) \geq f(z, x) + \langle w, y - x \rangle, \ \forall y \in C \}. $$
In particular,
$$ \partial_2 f(z , z)  = \{ w \in \mathbb{H} : f(z, y) \geq  \langle w, y - z \rangle, \ \forall y  \in C \}.$$
Let $\Omega$ be an open convex set containing $C$. The next lemma can be considered as an infinite-dimensional version of Theorem 24.5 in \cite{Roc}
\begin{lemma}\label{Lem2.2}\cite[Proposition 4.3]{VSH}
Let $f : \Omega \times \Omega \to \mathbb{R}$ be a function satisfying conditions $(A_1)$  on $\Omega$ and $(A_2)$ on $C$. Let $\bar{x}, \bar{y} \in \Omega$ and $\{x^k\}$, $\{y^k\}$ be two sequences in $\Omega$ converging weakly to $\bar{x}, \bar{y}$, respectively. Then, for any $\epsilon > 0$, there exist $\eta >0$ and $k_{\epsilon} \in \mathbb{N}$ such that
$$ \partial_2 f(x^k, y^k) \subset \partial_2 f(\bar{x}, \bar{y}) + \frac{\epsilon}{\eta}B, $$
for every $k \geq k_\epsilon$, where $B$ denotes the closed unit ball in $\mathbb{H}.$
\end{lemma}

\begin{lemma}\label{Lem2.6}  Suppose the bifunction $f$ satisfies the assumptions $(A_1)$ on $\Omega$ and $(A_2)$ on $C$. If $\{x^k \} \subset C $ is bounded, $\rho > 0$,  and $\{y^k\}$ is a sequence such that
$$ y^k = \arg\min\Big\{  f(x^k, y) +  \frac{\rho}{2} \| y-x^k\|^2: \ y\in C\Big\},$$
then $\{y^k\}$ is bounded.
  \end{lemma}
{\bf Proof.}
Firstly, we show that if $\{x^k \}$  converges weakly to $x^*$, then $\{y^k\}$ is bounded. Indeed,
$$ y^k = \text{arg}\min\Big\{  f(x^k, y) +  \frac{\rho}{2} \| y-x^k\|^2: \ y\in C \Big\},$$
and
$$ f(x^k, x^k) +  \frac{\rho}{2} \| y^k-x^k\|^2 = 0,$$
therefore
$$ f(x^k, y^k) +  \frac{\rho}{2} \| y^k-x^k\|^2 \leq 0, \ \forall k. $$
In addition, for all $w^k \in \partial_2 f(x^k, x^k)$ we have
 $$ f(x^k, y^k) +  \frac{\rho}{2} \| y^k-x^k\|^2  \geq \langle w^k, y^k - x^k \rangle + \frac{\rho}{2} \|y^k - x^k\|^2.$$
This implies $ -\|w^k\| \|y^k - x^k \| + \frac{\rho}{2} \|y^k - x^k\|^2 \leq 0$. Hence,
$$\|y^k - x^k\| \leq \frac{2}{\rho} \|w^k\|, \ \forall k.$$
 Because  $\{x^k\}$ converges weakly to $x^*$ and $w^k \in \partial_2f(x^k, x^k)$, by Lemma~\ref{Lem2.2}, the sequence $\{w^k\}$ is bounded, combining with the boundedness of $\{x^k\}$, we get $\{y^k\}$ is also bounded.

Now we prove the Lemma~\ref{Lem2.6}. Suppose that $\{y^k\}$ is unbounded, i.e., there exists an subsequence $\{y^{k_i}\} \subseteq  \{y^k\}$ such that $\lim_{i \to \infty}\|y^{k_i}\| = + \infty$. By the boundedness of  $\{x^k\}$, it implies  $\{x^{k_i}\}$ is also bounded, without loss of generality, we may assume that $\{x^{k_i}\}$ converges weakly to some $ x^*$. By the same argument as above, we obtain $\{y^{k_i}\}$ is bounded, which contradicts.
Therefore  $\{y^k\}$ is bounded.
\hfill$\Box$

\begin{lemma}\label{L2.5}\cite{KT}
Let $C$ be a nonempty closed convex subset of  $\mathbb{H}$ . Assume that $T$ is an
($\alpha, \beta, \gamma, \delta$)-symmetric generalized hybrid self-mapping of $C$ such that $Fix(T) \neq \emptyset$ and the conditions $(1) \ \alpha + 2\beta + \gamma \geq 0$,  $(2) \ \alpha + \beta > 0$ and $(3) \ \delta \geq 0 $ hold. Then T is quasi-nonexpansive.
\end{lemma}

\begin{lemma}\label{L2.7}\cite{HST}
Let $C$ be a nonempty closed convex subset of  $\mathbb{H}$ . Assume that $T$ is an
($\alpha, \beta, \gamma, \delta$)-symmetric generalized hybrid self-mapping of $C$ such that $Fix(T) \neq \emptyset$ and the conditions $(1) \ \alpha + 2\beta + \gamma \geq 0$,  $(2) \ \alpha + \beta > 0$ and $(3) \  \delta \geq 0 $ hold. Then $I - T$ is demiclosed at $0$, i.e., $x^k \rightharpoonup \bar{x}$ and $x^k - Tx^k \to 0$ imply $\bar{x} \in Fix(T)$.
\end{lemma}

\section{An extragradient algorithm}
 {\bf Algorithm 1}
\begin{itemize}
   \item[]{\bf Initialization.} Pick  $x^0 = x^g \in C$, choose parameters $\{ \rho_k \} \subset [ \b{$\rho$}, \  \bar{\rho} ]$, with $0 < \b{$\rho$} \leq \bar{\rho} < \min\{\frac{1}{2L_1}, \frac{1}{2L_2}\} $, \\
        \indent $ \{ \alpha_k \} \subset [0, 1], \ \lim_{k \to \infty} \alpha_k = 1$, $\{\beta_k\} \subset [0, \bar{\beta}] \subset [0, 1)$.
   \item[]{\bf Iteration $k$} (k = 0, 1, 2, ...). Having $x^k$  do the following steps:
	\begin{itemize}
	 \item[]{\it Step 1.} Solve the successively strongly convex programs
            $$\min\Big\{  \rho_k f(x^k, y) +  \frac{1}{2} \| y-x^k\|^2: \ y\in C\Big\} \ \eqno CP(x^k)$$
            $$ \min\Big\{  \rho_k f(y^k, y) +  \frac{1}{2} \| y-x^k\|^2: \ y\in C\Big\} \eqno CP(y^k, x^k)$$
             to obtain their unique solutions $y^k$ and $z^k$ respectively.
	 \item[]{\it Step 2.} Compute
                   \begin{align*}
                      t^{k} &= \alpha_k x^k + (1 -\alpha_k) Tx^k,\\ 
                     u^{k} &= \beta_k t^k + (1 -\beta_k) Tz^k. 
                      \end{align*}
	\item[]{\it Step 3.} Define \begin{align*}
                                          C_k & = \{x \in \mathbb{H}: \|x - u^k \| \leq \|x - x^k \| \},\\
                                          Q_k & =  \{x \in \mathbb{H}: \langle x - x^k , x^g  - x^k \rangle \leq 0 \}, \\
                                               A_k & = C_k \cap Q_k \cap C.
                                           \end{align*}
                                 \indent Take $x^{k+1} = P_{A_k}(x^g)$, and go to Step 1  with $k$ is replaced by  $k +1 $.\\
\end{itemize}
\end{itemize}

Before going to prove the convergence of  this algorithm, let us recall the following result which was proved in \cite{Anh}
\begin{lemma}\label{Lem3.1}\cite{Anh}
Suppose that $x^* \in \text{ Sol}(C, f)$,  then under assumptions $(A_2)$, $(A_3)$, and $(A_4)$, we have:
\begin{itemize}
\item[(i)]  $\rho_k[f(x^k, y)-f(x^k, y^k)] \geq \langle y^k-x^k, y^k-y \rangle,\;\forall y \in C.$
\item[(ii)]  $\| z^k-x^* \|^2  \leq \|x^k-x^* \|^2- (1 - 2\rho_k L_1) \| x^k-y^k \|^2 - (1-2 \rho_k L_2) \| y^k - z^k \|^2,\;\;\forall k.$
\end{itemize}
\end{lemma}

\begin{theorem}\label{The3.1} Suppose that the set $S = Sol(C, f) \cap Fix(T)$ is nonempty. Then  under assumptions ($A_1$), ($A_2$), $(A_3)$, $(A_4)$, and $(A_5)$ the sequences $\{x^k\}$, $\{y^k\}$, $\{z^k\}$ generated by Algorithm 1 converge strongly to the solution $x^* = P_S(x^g)$.
      \end{theorem}
{\bf Proof.}
Take  $q \in \ S$, i.e., $q \in  Sol(C, f) \cap Fix(T)$. By definition of $\rho_k$: $ 0 < \b{$\rho$} \leq \rho_k \leq{\bar{\rho}} < \text{min}\{\frac{1}{2L_1}, \frac{1}{2L_2}\}$, we get from Lemma~\ref{Lem3.1} that
\begin{equation}\label{3.3}
   \|z^k - q \| \leq \|x^k - q \|.
\end{equation}
By definition of $t^k$, we have
\begin{equation*}
\begin{aligned}
\|t^k - q\| & = \| \alpha_kx^k + (1 - \alpha_k) Tx^k - q \| \\
                & \leq \alpha_k\|x^k - q\| + (1 - \alpha_k)\|Tx^k - q \|.
\end{aligned}
\end{equation*}
Since $T$ is ($\alpha, \beta, \gamma, \delta)$-symmetric generalized hybrid mapping with $\alpha + 2\beta + \gamma \geq 0$, $\alpha + \beta > 0$, $\delta \geq 0$. From Lemma~\ref{L2.5} it is quasi-nonexpansive, so
\begin{equation}\label{3.4}
\| t^k - q \|  \leq  \| x^k - q \|.
\end{equation}
Similarly
\begin{equation}\notag
\begin{aligned}
\|u^k - q\| & = \|\beta_k t^k + (1 - \beta_k) Tz^k - q \|\\
&\leq \beta_k \|t^k - q \| + (1 - \beta_k) \|Tz^k - q \|\\
&\leq \beta_k \|x^k - q \| + (1 - \beta_k) \|z^k - q \|.
\end{aligned}
\end{equation}
Combining with (\ref{3.3}) yields
\begin{equation}\label{3.5}
\| u^k - q \|  \leq  \| x^k - q \|.
\end{equation}
Next, we show that $S \subset C_k \cap Q_k, \ \forall k.$ Indeed, from (\ref{3.5}) it implies that $q \in C_k$, or $S \subset C_k$ for all $k$. We prove $S \subset Q_k$ by induction. It is clear that $S \subset Q_0.$ If $S \subset Q_k$, i.e., $\langle q - x^k, x^g - x^k \rangle \leq 0, \ \forall q \in S.$ Since $x^{k+1} = P_{A_k}(x^g)$ we obtain $\langle x - x^{k+1}, x^g - x^{k+1} \rangle \leq 0, \ \forall x \in A_k .$ Especially, $\langle q - x^{k+1}, x^g - x^{k+1} \rangle \leq 0, \ \forall q \in S$. So $S \subset C_k \cap Q_k, \ \forall k.$ \\
From definition of $Q_k$, it implies that $x^k =  P_{Q_k}(x^g)$, so $\|x^k - x^g\| \leq \|x - x^g\|, \ \forall x \in Q_k$. In particular
\begin{equation}\label{3.6}
\|x^k - x^g\| \leq \|q - x^g \|, \ \forall k, \ \forall q \in S.
\end{equation}
Consequently, $\{x^k\}$ is bounded. Combining with (\ref{3.4}), (\ref{3.5}), we get $\{t^k\}$, $\{u^k\}$ are also bounded.\\
In addition,
\begin{equation}\notag
\begin{aligned}
\|x^{k+1} - x^k \|^2 & = \|x^{k+1} - x^g + x^g - x^k \|^2 \\
                                       & = \|x^{k+1} - x^g \|^2 + \|x^g - x^k \|^2 + 2 \langle x^{k+1} - x^g, x^g - x^k \rangle\\
                                       & = \|x^{k+1} - x^g \|^2 - \|x^g - x^k \|^2 + 2 \langle x^{k+1} - x^k, x^g - x^k \rangle.
\end{aligned}
\end{equation}
Since $x^{k+1} \in Q_k$, it implies from the above inequality that
\begin{equation}\label{3.7}
\|x^{k+1} - x^k \|^2 \leq \|x^{k+1} - x^g \|^2 - \|x^{k} - x^g \|^2.
\end{equation}
Therefore $\{\|x^k - x^g \|\}$ is nondecreasing sequence. In view of (\ref{3.6}), the limit \\
 $\lim_{k \to \infty}\|x^k - x^g\|$ exists. Hence, it also follows from (\ref{3.7}) that
\begin{equation}\label{3.8}
\lim_{k \to \infty}\|x^{k+1} - x^k \| = 0.
\end{equation}
Because $x^{k+1} \in C_k$, it implies that
\begin{equation}\notag
\begin{aligned}
\|u^{k} - x^k \| & \leq \|u^k - x^{k+1} \| + \|x^{k+1} - x^k \| \\
                               & \leq 2\|x^{k+1} - x^k \|,
\end{aligned}
\end{equation}
therefore, we deduce from (\ref{3.8}) that
\begin{equation}\label{3.9}
\lim_{k \to \infty} \|u^k - x^k \| = 0.
\end{equation}
Besides that $\lim_{k \to \infty} \alpha_k = 1$, so
\begin{equation}\label{3.10}
\lim_{k \to \infty} \|t^k - x^k \| = \lim_{k \to \infty} (1 - \alpha_k)\|x^k - T x^k \| = 0.
\end{equation}
It is clear that
\begin{equation}\notag
\begin{aligned}
\| u^k - q \|^2 & = \| \beta_k (t^k - q) + (1 - \beta_k) (Tz^k - q) \|^2\\
                           & =   \beta_k \|t^k - q\|^2 + (1 - \beta_k) \|Tz^k - q\|^2 -  \beta_k(1 - \beta_k) \| t^k - Tz^k \|^2 \\
                            & \leq  \beta_k \|t^k - q\|^2 + (1 - \beta_k) \|Tz^k - q\|^2.
\end{aligned}
\end{equation}
In view of (\ref{3.4}), Lemma~\ref{L2.5}, and Lemma~\ref{Lem3.1}, yields
\begin{equation*}
\| u^k - q \|^2 \leq   \|x^k - q\|^2 - (1 - \beta_k) \big[ (1 - 2\rho_k L_1) \| x^k-y^k \|^2 - (1-2 \rho_k L_2) \| y^k - z^k \|^2  \big].
\end{equation*}
Hence
\begin{equation}\label{3.11}
(1 - \beta_k)\big[(1 - 2\rho_k L_1) \| x^k-y^k \|^2 + (1-2 \rho_k L_2) \| y^k - z^k \|^2 \big]  \leq \|x^k - u^k \| \big[ \|x^k - q \| + \| u^k - q \| \big]
\end{equation}
Since $0 < 1 -\bar{\beta} \leq 1- \beta_k$; $0 < \b{$\rho$} \leq \rho_k \leq \bar{\rho} < \min\{\frac{1}{2L_1}, \frac{1}{2L_2}\}$, and (\ref{3.9}), we get from (\ref{3.11}) that
\begin{equation}\label{3.12}
\lim_{k \to \infty} \|x^k - y^k \| = 0.
\end{equation}
\begin{equation}\label{3.13}
\lim_{k \to \infty} \|y^k - z^k \| = 0.
\end{equation}
\begin{equation}\label{3.14}
\lim_{k \to \infty} \|x^k - z^k \| = 0.
\end{equation}
By definition of $u^k$, we have $(1 - \beta_k) Tz^k =  u^k - \beta_k t^k$. Hence
\begin{equation*}
\begin{aligned}
(1 - \bar{\beta}) \|Tz^k - z^k \| & \leq \|(1- \beta_k) Tz^k - (1 - \beta_k)z^k \| \\
                                                         & = \|u^k - z^k - \beta_k(t^k - z^k)\| \\
                                                         & \leq \| u^k - z^k \| + \beta_k\|t^k - z^k\| \\
                                                         & \leq \| u^k - x^k \| + \beta_k\|t^k - x^k\|  + (1 + \beta_k) \|x^k - z^k\|.
\end{aligned}
\end{equation*}
Combining this fact with (\ref{3.9}), (\ref{3.10}), and  (\ref{3.14}) we receive in the limit that
\begin{equation}\label{3.15}
\lim_{k\to \infty} \|Tz^k - z^k \| = 0.
\end{equation}

Next we show that any weak accumulation point of $\{x^k\}$ belongs to $S$. Indeed, suppose that $\{x^{k_i}\} \subset \{x^k\}$ and $x^{k_i} \rightharpoonup p$  as $i \to \infty$. From (\ref{3.12}), (\ref{3.13}), and (\ref{3.14}) we get $y^{k_i} \rightharpoonup p$, and $z^{k_i} \rightharpoonup p$ as $i \to \infty$.

Replacing $k$ by $k_i$ in assertion (i) of Lemma~\ref{Lem3.1} we get
\begin{equation*}
\rho_{k_i}\big[ f(x^{k_i}, y) - f(x^{k_i}, y^{k_i}) \big] \geq  \langle x^{k_i} - y^{k_i} , y - y^{k_i} \rangle, \ \forall y \in C.
\end{equation*}
Hence
\begin{equation}\label{3.16}
\rho_{k_i}\big[ f(x^{k_i}, y) - f(x^{k_i}, y^{k_i}) \big] \geq   -\|x^{k_i} - y^{k_i}\| \| y - y^{k_i}\|.
\end{equation}
Letting $i \to \infty$, by jointly weak continuity of $f$ and  (\ref{3.12}), we obtain in the limit from (\ref{3.16}) that
$$ f(p, y) - f(p, p) \geq 0. $$
So
$$ f(p,y)  \geq 0, \ \forall y \in C, $$
which  means that $p$ is a solution of EP($C, f$).\\
By (\ref{3.15}), we have that $\lim_{i \to \infty}\|Tz^{k_i} - z^{k_i}\| = 0$. Since $z^{k_i} \rightharpoonup p$ and demiclosedness at zero of $I - T$, Lemma~\ref{L2.7}, we get $Tp = p$, i.e., $p \in Fix(T).$ \\
Hence $p \in S$.

Now, we set $x^* = P_{S}(x^g)$. From (\ref{3.6}) one has,
\begin{equation*}\label{4.6}
\|x^k - x^g\| \leq \|x^* - x^g\|, \ \forall k.
\end{equation*}
It is immediate from Lemma~\ref{Lem2.4} that $x^k$ converges strongly to $x^*$. Combining with (\ref{3.12}), (\ref{3.14}) we have that $y^k$, $z^k$ converge strongly to $x^*.$ This completes the proof.
\hfill$\Box$


\section{An extragradient algorithm with linesearch}
{\bf Algorithm 2}
\begin{itemize}
\item[]{\bf Initialization.} Pick  $x^0 = x^g \in C$, choose parameters $\eta, \mu \in (0, 1); \ 0 < \b{$\rho$} \leq \bar{\rho}$, $\{ \rho_k \} \subset [ \b{$\rho$}, \ \bar{\rho}]$; \\
         \indent $ \{ \alpha_k \} \subset [0, 1],$  $\lim_{k \to \infty} \alpha_k = 1$;    $\{\beta_k\} \subset [0, \bar{\beta}] \subset [0, 1)$; $\gamma_k \in [\b{$\gamma$}, \bar{\gamma}] \subset (0, 2)$.
\item[]{\bf Iteration $k$} (k = 0, 1, 2, ...). Having $x^k$  do the following steps:
      \begin{itemize}
		\item[]{\it Step 1.} Solve the strongly convex program
		$$\min\Big\{  \rho_k f(x^k, y) +  \frac{1}{2} \| y-x^k\|^2: \ y\in C\Big\}  \ \ \  \ \ \ \ \ \ \ \ \      	\ \eqno CP(x^k)$$
			to obtain its unique solutions $y^k$. \\
			If $y^k = x^k$, then set $v^k = x^k$. Otherwise go to Step 2.
	\item[]{\it Step 2.} (Armijo linesearch rule) Find $m_k$ as the smallest positive integer 	   number $m$ such that
	\begin{equation}\label{4.1}
	\begin{cases}
	z^{k,m} = (1 - \eta^m)x^k + \eta^my^m \\
	f(z^{k,m}, x^k) - f(z^{k,m}, y^k)  \geq \frac{\mu}{2 \rho_k}\|x^k - y^k\|^2.
	\end{cases}
	\end{equation}
	Set $\eta_k = \eta^{m_k}$, $z^k = z^{k, m_k}$.
	\item[]{\it Step 3.} Select $w^k \in \partial_2f(z^k, x^k)$, and compute $v^k = P_C(x^k - \gamma_k.\sigma_k.w^k)$, \\
     where $\sigma_k = \frac{f(z^k, x^k)}{\|w^k\|^2}$.
     \item[]{\it Step 4.} Compute
       \begin{align*}
	t^{k} &= \alpha_k x^k + (1 -\alpha_k) Tx^k,\\ 
	u^{k} &= \beta_k t^k + (1 -\beta_k) Tv^k. 
          \end{align*}
    \item[]{\it Step 5.} Define \begin{align*}
	                          C_k & = \{x \in \mathbb{H}: \|x - u^k \| \leq \|x - x^k \| \},\\
	                          Q_k & =  \{x \in \mathbb{H}: \langle x - x^k , x^g  - x^k \rangle \leq 0 \}, \\
                               A_k & = C_k \cap Q_k \cap C.
                          \end{align*}
\indent Take $x^{k+1} = P_{A_k}(x^g)$, and go to Step 1  with $k$ is replaced by  $k +1 $.\\
\end{itemize}
\end{itemize}
 \begin{remark}\label{Rem41}
\begin{itemize}
\item[(i)] If  $y^k = x^k$  then $x^k$ is a solution to \text{EP}($C, f$);
\item[(ii)] If $y^k = x^k = t^k$ and $\alpha_k < 1$ or $y^k = x^k = u^k$, then $x^k \in Sol(C, f) \cap Fix(T)$.
\end{itemize}
\end{remark}
Firstly,  let us recall the following lemma which was proved in \cite{MQH}
\begin{lemma}\label{Lem4.1}\cite{MQH}
Suppose that $p \in \text{ Sol}(C, f)$, then under assumptions $(A_2)$, $(A_3)$, and $(A_4)$, we have:
\begin{itemize}
\item[(a)]  The linesearch is well defined;
\item[(b)]  $ f(z^k, x^k) > 0$;
\item[(c)] $0 \not\in \partial_2f(z^k, x^k)$;
\item[(d)] \begin{equation}\label{4.2} \|v^k - p\| \leq \|x^k - p\|^2 - \gamma_k( 2 - \gamma_k)(\sigma_k\|w^k\|)^2.
                 \end{equation}
\end{itemize}
\end{lemma}

\begin{theorem}\label{The4.1} Suppose that the set $ S = Sol(C, f) \cap Fix(T)$ is nonempty, the bifunction $f$ satisfies assumptions ($A_1$) on $\Omega$, ($A_2$), and $(A_3)$ on $C$, the mapping $T$ satisfies assumption $(A_5)$. Then the sequences $\{x^k\}$, $\{u^k\}$ generated by Algorithm 2 converge strongly to the solution $x^* = P_S(x^g)$.
      \end{theorem}


{\bf Proof.}
Take $q \in S$. Since $\gamma_k \in [\b{$\gamma$}, \bar{\gamma}] \subset (0, 2)$, we get from Lemma~\ref{Lem4.1} that
\begin{equation}\label{4.3}
   \|v^k - q \| \leq \|x^k - q \|.
\end{equation}
By definition of $t^k$, we have
\begin{equation*}
\begin{aligned}
\|t^k - q\|^2 & = \| \alpha_k(x^k - q) + (1 - \alpha_k) (Tx^k - q) \|^2 \\
                & = \alpha_k\|x^k - q\|^2 + (1 - \alpha_k)\|Tx^k - q \|^2 - \alpha_k(1 - \alpha_k)\|Tx^k - x^k \|^2.
\end{aligned}
\end{equation*}
Since $T$ is a ($\alpha, \beta, \gamma, \delta)$-symmetric generalized hybrid mapping with $\alpha + 2\beta + \gamma \geq 0$, $\alpha + \beta > 0$, $\delta \geq 0$. By Lemma~\ref{L2.5} it is quasi-nonexpansive, so
\begin{equation}\label{4.4}
\| t^k - q \|  \leq  \| x^k - q \|.
\end{equation}
Similarly,
\begin{equation}\label{4.5}
\| u^k - q \|  \leq  \| x^k - q \|.
\end{equation}
Next, we show that $S \subset C_k \cap Q_k, \ \forall k.$ Indeed, from (\ref{4.5}) it implies that $q \in C_k$, or $S \subset C_k$. We prove $S \subset Q_k$ by induction, it is clear that $S \subset Q_0.$ If $S \subset Q_k$, i.e., $\langle q - x^k, x^g - x^k \rangle \leq 0, \ \forall q \in S.$ Since $x^{k+1} = P_{A_k}(x^g)$ we obtain $\langle x - x^{k+1}, x^g - x^{k+1} \rangle \leq 0, \ \forall x \in A_k .$ Especially, $\langle q - x^{k+1}, x^g - x^{k+1} \rangle \leq 0, \ \forall q \in S$. So $S \subset C_k \cap Q_k, \ \forall k.$ \\
From definition of $Q_k$, it implies that $x^k =  P_{Q_k}(x^g)$, so $\|x^k - x^g\| \leq \|x - x^g\|, \ \forall x \in Q_k$. In particular
\begin{equation}\label{4.6}
\|x^k - x^g\| \leq \|q - x^g \|, \ \forall k, \ \forall q \in S.
\end{equation}
Consequently, $\{x^k\}$ is bounded. Combining with (\ref{4.4}), (\ref{4.5}), we get $\{t^k\}$, $\{u^k\}$ are also bounded.\\
In addition,
\begin{equation}\notag
\begin{aligned}
\|x^{k+1} - x^k \|^2 & = \|x^{k+1} - x^g + x^g - x^k \|^2 \\
                                       & = \|x^{k+1} - x^g \|^2 + \|x^g - x^k \|^2 + 2 \langle x^{k+1} - x^g, x^g - x^k \rangle\\
                                       & = \|x^{k+1} - x^g \|^2 - \|x^g - x^k \|^2 + 2 \langle x^{k+1} - x^k, x^g - x^k \rangle.
\end{aligned}
\end{equation}
Since $x^{k+1} \in Q_k$, it implies from the above inequality that
\begin{equation}\label{4.7}
\|x^{k+1} - x^k \|^2 \leq \|x^{k+1} - x^g \|^2 - \|x^{k} - x^g \|^2.
\end{equation}
Therefore $\{\|x^k - x^g \|\}$ is nondecreasing sequence. Together with (\ref{4.6}), the limit
 $\lim_{k \to \infty}\|x^k - x^g\|$ does exist.\\
  Hence, it also follows from (\ref{4.7}) that
\begin{equation}\label{4.8}
\lim_{k \to \infty}\|x^{k+1} - x^k \| = 0.
\end{equation}
Because $x^{k+1} \in C_k$, it implies that
\begin{equation}\notag
\begin{aligned}
\|u^{k} - x^k \| & \leq \|u^k - x^{k+1} \| + \|x^{k+1} - x^k \| \\
                               & \leq 2\|x^{k+1} - x^k \|,
\end{aligned}
\end{equation}
therefore, we deduce from (\ref{4.8}) that
\begin{equation}\label{4.9}
\lim_{k \to \infty} \|u^k - x^k \| = 0.
\end{equation}
Besides that $\lim_{k \to \infty} \alpha_k = 1$, so
\begin{equation}\label{4.10}
\lim_{k \to \infty} \|t^k - x^k \| = \lim_{k \to \infty} (1 - \alpha_k)\|x^k - T x^k \| = 0.
\end{equation}
It is clear that
\begin{equation}\notag
\begin{aligned}
\| u^k - q \|^2 & = \| \beta_k (t^k - q) + (1 - \beta_k) (Tv^k - q) \|^2\\
                           & =   \beta_k \|t^k - q\|^2 + (1 - \beta_k) \|Tv^k - q\|^2 -  \beta_k(1 - \beta_k) \| t^k - Tv^k \|^2 \\
                            & \leq  \beta_k \|t^k - q\|^2 + (1 - \beta_k) \|Tv^k - q\|^2.
\end{aligned}
\end{equation}
In view of (\ref{4.4})  and Lemma~\ref{Lem4.1}, yields
\begin{equation*}
\| u^k - q \|^2 \leq   \|x^k - q\|^2 - (1 - \beta_k) \gamma_k(2 - \gamma_k)(\sigma_k \|w^k\|)^2.
\end{equation*}
Hence
\begin{equation}\label{4.11}
(1 - \beta_k) \gamma_k(2 - \gamma_k)(\sigma_k \|w^k\|)^2  \leq \|x^k - u^k \| \big[ \|x^k - q \| + \| u^k - q \| \big].
\end{equation}
Since $0 < 1 -\bar{\beta} \leq 1- \beta_k$; $\gamma_k \in [\b{$\gamma$}, \bar{\gamma}] \subset (0, 2)$, and (\ref{4.9}), we get from (\ref{4.11}) that
\begin{equation}\label{4.12}
\lim_{k \to \infty} \sigma_k\|w^k \| = 0.
\end{equation}
Because $v^k = P_C(x^k - \gamma_k\sigma_k w^k)$, one has
$$ \|v^k - x^k\| \leq  \gamma_k\sigma_k \|w^k\|.$$
Combining with (\ref{4.12}) we get
\begin{equation}\label{4.14}
\lim_{k \to \infty} \|v^k - x^k \| = 0.
\end{equation}
By definition of $u^k$, we have $(1 - \beta_k) Tv^k =  u^k - \beta_k t^k$.
Hence
\begin{equation*}
\begin{aligned}
(1 - \bar{\beta}) \|Tv^k - v^k \| & \leq \|(1- \beta_k) Tv^k - (1 - \beta_k)v^k \| \\
                                                         & = \|u^k - v^k - \beta_k(t^k - v^k)\| \\
                                                         & \leq \| u^k - v^k \| + \beta_k\|t^k - v^k\| \\
                                                         & \leq \| u^k - x^k \| + \beta_k\|t^k - x^k\|  + (1 + \beta_k) \|x^k - v^k\|.
\end{aligned}
\end{equation*}
Combining this fact with (\ref{4.9}), (\ref{4.10}), (\ref{4.14}), we receive in the limit that
\begin{equation}\label{4.15}
\lim_{k\to \infty} \|Tv^k - v^k \| = 0.
\end{equation}
By Lemma~\ref{Lem2.6}, $\{y^k\}$ is bounded, consequently $\{z^k\}$ is bounded. From Lemma~\ref{Lem2.2}, $\{w^k\}$ is bounded.
In view of  (~\ref{4.12}) yields
\begin{equation}\label{4.16}
  \lim_{k \to \infty} f(z^k, x^k) = \lim_{k \to \infty} [\sigma_k\|w^k\|] \|w^k\| = 0.
\end{equation}
We have
\begin{equation*}
\begin{aligned}
0 = f(z^k, z^k) & = f(z^k, (1 - \eta_{k}) x^k + \eta_k y^k) \\
                        & \leq (1 - \eta_k)f(z^k, x^k) + \eta_k f(z^k, y^k),
\end{aligned}
\end{equation*}
so, we get from (\ref{4.1}) that

\begin{equation*}
\begin{aligned}
f(z^k, x^k) & \geq \eta_k [f(z^k, x^k) - f(z^k, y^k)]\\
                  &\geq  \frac{\mu}{2\rho_k} \eta_k \|x^k - y^k \|^2.
\end{aligned}
\end{equation*}
Combining with (\ref{4.16}) one has
\begin{equation}\label{4.17}
 \lim_{k \to \infty} \eta_k\| x^k - y^k\|^2 = 0.
\end{equation}
Next, we show that any weak accumulation point of $\{x^k\}$ belongs to $S$. Indeed, suppose that $\{x^{k_i}\} \subset \{x^k\}$ and $x^{k_i} \rightharpoonup p$ as $i \to \infty$. \\
From (\ref{4.17}) we get
\begin{equation}\label{4.18}
 \lim_{i \to \infty} \eta_{k_i}\| x^{k_i} - y^{k_i}\|^2 = 0.
\end{equation}
We now consider two distinct cases:

{\it Case 1.} $\lim\sup_{i \to \infty}\eta_{k_i} > 0$. \\
 Without loss of generality, we may assume that there exists $\bar{\eta} > 0$ such that  $ \eta_{{k}_i} > \bar{\eta }, \ \forall i \geq i_0 $,  use this fact  and  from (\ref{4.18}),  one has
\begin{equation}\label{4.19}
\lim_{i \to \infty}{\|x^{{k}_i}-y^{{k}_i}\|} = 0.
\end{equation}
 Remember that  $ x^k \rightharpoonup p$, together with (\ref{4.19}), it implies  that $ y^{k_i} \rightharpoonup p$ as $i \to \infty$. \\

From assertation (i) of Lemma~\ref{Lem3.1} we get
\begin{equation*}
\rho_{k_i}\big[ f(x^{k_i}, y) - f(x^{k_i}, y^{k_i}) \big] \geq  \langle x^{k_i} - y^{k_i} , y - y^{k_i} \rangle, \ \forall y \in C.
\end{equation*}
Hence
\begin{equation}\label{4.20}
\rho_{k_i}\big[ f(x^{k_i}, y) - f(x^{k_i}, y^{k_i}) \big] \geq   -\|x^{k_i} - y^{k_i}\| \| y - y^{k_i}\|.
\end{equation}
Letting $i \to \infty$, by jointly weak continuity of $f$ and  (\ref{4.19}), we obtain in the limit from (~\ref{420}) that
$$ f(p, y) - f(p, p) \geq 0. $$
So
$$ f(p,y)  \geq 0, \ \forall y \in C, $$
which  means that $p$ is a solution of EP($C, f).$

{\it Case 2.} \   $\lim_{i \to \infty}{\eta_{k_i}} = 0$.\\
From the boundedness of $\{y^{k_i}\}$, without loss of generality we may assume that $ y^{k_i} \rightharpoonup  \bar{y}$ as $i \to \infty$. \\
 Replacing $y$ by $x^{k_i}$ in ($i$) of Lemma~\ref{Lem3.1} we get
\begin{equation}\label{4.21}
f(x^{k_i}, y^{k_i}) \leq -\frac{1}{\rho_{k_i}} \| y^{k_i} - x^{k_i} \|^2.
\end{equation}
In the other hand, by the Armijo linesearch rule (\ref{4.1}), for $m_{k_i} - 1$, we have
\begin{equation}\label{4.22}
 f(z^{k_i, m_{k_i} - 1}, x^{k_i}) - f(z^{k_i, m_{k_i} - 1}, y^{k_i})  < \frac{\mu}{2\rho_{k_i}} \| y^{k_i}-x^{k_i}\|^2.
\end{equation}
Combining with (\ref{4.21}) we get
\begin{equation}\label{4.23}
f(x^{k_i}, y^{k_i}) \leq -\frac{1}{\rho_{k_i}} \| y^{k_i} - x^{k_i} \|^2 \leq \frac{2}{\mu} \big[f(z^{k_i, m_{k_i} - 1}, y^{k_i}) - f(z^{k_i, m_{k_i} - 1}, x^{k_i}) \big].
\end{equation}
 According to the algorithm, we have $z^{k_i, m_{k_i} - 1} = (1-\eta^{m_{k_i} - 1})x^{k_i} + \eta^{m_{k_i} - 1}y^{k_i}$, $\eta^{k_i, m_{k_i} - 1} \to 0$  and $x^{k_i} $ converges weakly to $p$, $y^{k_i}$ converges weakly to $\bar{y}$, it implies that $z^{k_i, m_{k_i} - 1} \rightharpoonup p$ as $i \to \infty$. Beside that $\{\frac{1}{\rho_{k_i}}\|y^{k_i} - x^{k_i}\|^2\}$ is bounded, without loss of generality, we may assume that $\lim_{i \to +\infty}\frac{1}{\rho_{k_i}}\|y^{k_i} - x^{k_i}\|^2$ exists. Hence, we get in the limit from (\ref{4.23}) that
\begin{equation*}
f(p, \bar{y}) \leq -  \lim_{i \to +\infty}\frac{1}{\rho_{k_i}}\|y^{k_i} - x^{k_i}\|^2 \leq \frac{2}{\mu}f(p, \bar{y}).
\end{equation*}
Therefore, $f(p, \bar{y}) = 0$ and $\lim_{i \to +\infty}\|y^{k_i} - x^{k_i}\|^2 = 0$. By the Case 1,
it is immediate that $p$ is a solution of EP($C, f$).\\
In addition, from (\ref{4.14}) and (\ref{4.15}), we have $v^{k_i} \rightharpoonup p$ and $\lim_{i \to \infty}\|Tv^{k_i} - v^{k_i}\| = 0$. By Lemma~\ref{L2.7}, $I - T$ is demiclosed at zero, we get $Tp = p$, i.e., $p \in Fix(T)$.\\
Hence $p \in S$.\\
Now, we set $x^* = P_{S}(x^g)$. From (\ref{4.6}) one has,
\begin{equation}\label{4.25}
\|x^k - x^g\| \leq \|x^* - x^g\|, \ \forall k.
\end{equation}
We get from Lemma~\ref{Lem2.4} that $x^k$ converges strongly to $x^*$. Combining with (\ref{4.9}) we  also have that $u^k$ converges strongly to $x^*.$ The proof is completed.
\hfill$\Box$

 {\bf Conclusion}. We have introduced two iterative methods for finding a common point of the solution set of a pseudomonotone equilibrium problem and the set of fixed points of a symmetric generalized hybrid mapping in a real Hilbert space. The basic iteration used in this paper is the extragradient  iteration with or without the incorporation of a linesearch procedure. The strong convergence of the iterates has been obtained.

\begin{acknowledgements}
This research was supported by Basic Science Research Program through the National Research Foundation of Korea (NRF) funded by the Ministry of Education (NRF-2013R1A1A2A10008908). The first author is supported in part by NAFOSTED, under the project 101.01-2014-24.
\end{acknowledgements}


\end{document}